\documentclass[12pt,emtex]{article}
\usepackage{amsfonts}
\usepackage{amssymb}
\usepackage{amscd}
\usepackage{amsmath}
\usepackage{amsthm}

\newcommand{\BEQ}{\begin{equation}}
\newcommand{\EEQ}{\end{equation}}

\newtheorem{df}{Definition}
\newtheorem{Th}{Theorem}
\newtheorem{prop}{Proposition}

\newtheorem{lemma}{Lemma}
\newtheorem{cor}{Corollary}

\newtheorem{rem}{Remark}

\newtheorem{Rem}{Remark}

\def\nn{\nonumber}

\def\bea{\begin{eqnarray}}
\def\eea{\end{eqnarray}}

\def\o{\omega}
\def\w{\wedge}

\def\a{\alpha}
\def\b{\beta}
\def\g{\gamma}
\def\po{\pi^{\#}}

\newcommand{\eqdef}{\stackrel{\rm def}{=}}
\newcommand{\lc}{\mathcal L}
\newcommand{\mc}{\mathcal M}
\newcommand{\tlc}{\tilde{\mathcal L}}
\newcommand{\tim}{\tilde\imath}
\newcommand{\tid}{\tilde d}
\newcommand{\somc}{\Sigma\Omega(\mc)}

\newcommand{\ta}{\tilde\alpha}
\newcommand{\tb}{\tilde\beta}
\newcommand{\tg}{\tilde\gamma}
\newcommand{\td}{\tilde\delta}

\begin{document}

\begin{titlepage}
\hfill ITEP-TH-28/04 \vskip 2.5cm

\centerline{\LARGE \bf On the Lie-formality of Poisson manifolds}

\vskip 1.0cm \centerline{G.Sharygin \footnote{E-mail:
sharygin@itep.ru},  D.Talalaev \footnote{E-mail:
talalaev@itep.ru}}

\centerline{\sf Institute of Theoretical and Experimental Physics
\footnote{ITEP, 25 ul. B. Cheremushkinskaya, Moscow, 117259 Russia}}

\vskip 2.0cm

\centerline{\large \bf Abstract} Starting from the problem of describing cohomological invariants of Poisson manifolds
we prove in a sense a ``no-go'' result: the differential graded Lie algebra
of de Rham forms on a smooth Poisson manifold is formal.
\vskip 1.0cm

\end{titlepage}


\section{Introduction}
The notion of formality was first introduced in the context of
differential graded algebras in the beginning of 1970-ies. It was related to
the investigations in the rational homotopy theory. In brief, to say that a
DG algebra $A$ is formal would mean that $A$ is homotopy equivalent (as algebra)
to its cohomology $H(A)$ (with zero differential), i.e. that there exists a
sequence of quasi-isomorphisms of DG algebras, beginning at $A$ and
ending at $H(A)$:
\bea
\label{eq1}
\begin{array}{ccccccccccccc}
& &A_1& & & &A_3& & & & A_n& & \\
&\swarrow& &\searrow& &\swarrow& &\searrow& &\swarrow & &\searrow & \\
A=A_0& & & & A_2 & & & & \ldots & & & & A_{n+1}=H(A)
\end{array}
\eea


Here quasi-isomorphism is a homomorphism of DG algebras, inducing an
isomorphism of their cohomology.

More generally, formality is a particular case of the notion of
homotopy equivalence of DG algebras: one says that algebras $A$ and $B$ are
homotopy equivalent, if there exists a sequence of algebras and
homomorphisms, similar to \eqref{eq1}, connecting them. All homotopical constructions,
applied to homotopy equivalent algebras give equivalent results. In
particular, their (usual) cohomology are isomorphic as well as the
DG-algebraic cohomology with coefficients in a suitable module.

However, the contrary is not true in general: two algebras with isomorphic
cohomology are not necessarily homotopy equivalent. One should also note, that there can be no quasi-isomorphism $A\to B$ even if $A$
and $B$ are homotopy equivalent:
it is not always possible to find homotopy inverse homomorphisms of DG-algebras for all the arrows, pointing to the left in the sequence \eqref{eq1}. (Recall, that two chain maps $f:A\to B$ and $g:B\to A$ are called homotopy inverse to each other,
if both their compositions $fg$ and $gf$ are homotopic to the corresponding identity maps.) On the contrary, it is
always possible to find a homotopy inverse chain map for a quasi-isomorphism of chain complexes (i.e for a map inducing an isomorphism of cohomology), if the
characteristic of the ground field is zero.

Using the last observation, one can define a suitable generalization of the morphism of
DG-algebras, that would be homotopy invertible whenever it establishes an isomorphism
of cohomology (in the characteristic zero case). It is the strong
homotopical morphisms of algebras. In brief, such a morphism from $A$ to $B$ is a collection
of linear maps $F=\{f_i\}_{i\ge 1},\ f_i:A^{\otimes i}\to B$, where every
map $f_n$ is a chain homotopy between zero and a certain combination of the
maps $f_i$ with smaller indices. The first two equations of this series are
\begin{align}
\label{eq2}
d_Bf_1(a)&=f_1(d_Aa),\\
\label{eq3}
d_Bf_2(a,b)+f_2(d_Aa,b)-(-1)^{|a|}&f_2(a,d_Ab)=f_1(a\/b)-f_1(a)f_1(b),
\end{align}
where $d_A,\ d_B$\ are the differentials in $A,\ B$, and $a,b\in A$\ are
the arbitrary elements. These equations mean that $f_1$ is a
chain map, and that $f_2$ is a chain homotopy, making $f_1^*:H(A)\to H(B)$
into a homomorphism of algebras.

In particular, every homomorphism $f$ of
algebras can be represented as an $A_\infty$ morphism, by choosing $f_i=0,\
i\ge2$.

As it was mentioned above, an important property of $A_\infty$ maps is that they can be
inverted up to a homotopy
if $f_1$ establishes an isomorphism in cohomology. For instance,
every quasi-isomorphism of algebras can be inverted in the class of such
maps. Starting from this observation one can prove that two DG-algebras $A$ and $B$ are
homotopy-equivalent iff there exists an $A_\infty$-morphism $F:A\to B$
with quasi-isomorphic $f_1$. In this case, there will exist a homotopy inverse morphism
$G=F^{-1}:B\to A$. In particular, $A$ is formal, iff there exists an
$A_\infty$-quasi-isomorphism $F:A\to H(A)$ (or $G:H(A)\to A$).

In analogy with associative algebras, there exists homotopy theory of
differential graded Lie algebras (DGL-algebras). The analogue of $A_\infty$
morphisms in this case is called $L_\infty$. An $L_\infty$-map $F:A\to B$
(where $A$ and $B$ are Lie algebras) consists of a series of linear mappings
$f_i:\wedge^iA\to B$ (where $\wedge^nA$ denotes the $n$\/-th
external power of the graded space $A$). These mappings are subject to some
relations, analogous to the above (see \S3 below). In particular, one has the following
analog of \eqref{eq3}
\begin{equation}
\label{eq4}
d_Bf_2(a,b)-f_2(d_Aa,b)+(-1)^{|a|}f_2(a,d_Ab)=f_1(\{a,b\}_A)-\{f_1(a),f_1(b)\}_B,
\end{equation}
where $\{,\}_A,\ \{,\}_B$\ are the commutators in $A$ and $B$. One can prove
(see e.g. \cite{Kon1}), that, similarly to the DG-algebras case, any
$L_\infty$ map with a quasi-isomorphic first stage, can be homotopy inverted in the
class of $L_\infty$ morphisms. In particular, a DGL-algebra is formal, iff
there exists an $L_\infty$ morphism from its cohomology Lie algebra to
itself, with quasi-isomorphic $f_1$.

As in the case of associative algebras, homotopy equivalence of two DGLa's
implies homotopy equivalence of functorial constructions, that one applies
to them. For example, their DGL cohomology with coefficients in a (differential graded) module
are isomorphic; e.g. for a formal Lie algebra $A$, one can calculate its Lie
algebraic cohomology by substituting $H(A)$ (with trivial differential) for
$A$ in the standard complex. This makes such cohomology calculable, since,
for an arbitrary DG Lie algebra $A$, $H(A)$ is often of finite type (i.e.
has finite dimensional homogeneous components).

In this note we prove the formality of the DG Lie algebra of de Rham forms
on a Poisson manifold. This fact implies, that there are no cohomological
invariants, associated to this algebra, except for its (usual) cohomology.
Speaking somewhat loosely, this Lie algebra knows not more about the Poisson
structure on the manifold, than its cohomology does.

The composition of this note is as follows. In the next section we
recall the definition and the main properties of the DG Lie algebra of the
differential forms on a Poisson manifold (the so-called Brylinski complex).
We prove that the induced Lie bracket on the cohomology is trivial and
deduce the formula \eqref{eq5 4/5}, which is crucial for our
proof of the formality.

In \S3 we discuss the definitions and the basic properties of the $L_\infty$ algebras and
$L_\infty$ morphisms. Finally, in \S4 we prove the formality result for
Brylinski complex of a Poisson manifold.

{\bf Acknowledgements.}
We are indebted to A. Khoroshkin for the discussions,
and we are thankful to I. Nikonov for pointing
the formula \eqref{eq5 4/5}.
The  work of both authors has been partially supported by the RFBR grant
04-01-00702.

\section{Brylinski complex of a Poisson manifold}
Let us consider a Poisson manifold $(M,\pi).$ Here $\pi$ is a Poisson
bivector, which in particular implies that the bracket on $C(M)$
defined by the formula
$$\{f,g\}=\pi(df,dg)$$
endows this space with the Lie algebra structure. That is, this bracket is anti-commutative and the
Jacobi identity holds:
$$\{f,\{g,h\}\}=\{\{f,g\},h\}+\{g,\{f,h\}\}.$$
The latter equation is equivalent to the condition that
the Schouten square of $\pi$ vanishes:
$$\{\pi,\pi\}=0.$$
{\Rem\rm The Schouten bracket on polyvector fields is a unique graded
bracket which coincides with the commutator of vector fields being
restricted to the first graded component and obeying the Leibnitz rule.}

In this section we collect several known and rather simple facts
about the DGLA structure on the complex of de Rham forms on a Poisson
manifold. We give it here mostly for the methodological reasons.
Let
$$0\longrightarrow C(M)\stackrel{d}{\longrightarrow}
\Omega^1(M)\stackrel{d}{\longrightarrow}\Omega^2(M)\cdots $$
be the de Rham
complex of $M.$
The material of this section is organized as follows

\begin{itemize}
  \item Firstly we recall the construction of the Koszul bracket on the space of $1$-forms on a
Poisson manifold and prove that it endows $\Omega^1(M)$ with the lie algebra
structure.
  \item  Further we extend this bracket to the exterior algebra $\bigwedge^*
\Omega^1(M)$ and prove that it descends correctly to the de Rham complex.
  \item We prove that this generalized bracket behaves well with respect to the external differential and
that the induced Lie algebra structure on the cohomology is trivial (i.e. the
Lie bracket is equal to zero).
\end{itemize}

\paragraph{\fbox{Koszul bracket}}
on $1$-forms is given by the formula:
\bea \label{p}
\{\a,\b\}=-d\pi(\a,\b)+\lc_{\po\a}\b-\lc_{\po\b}\a
\eea
where $\lc$ denotes the Lie derivative.
On exact forms $df,dg$ it gives
$$\{df,dg\}=d\{f,g\}$$ where on the right-hand side we use the standard
Poisson bracket of functions on Poisson manifold.

Let us recall that there is a Lie algebroid structure on $T^*(M)$ with the
anchor map
$$\po:T^*(M)\rightarrow T(M).$$
The important property of a Lie
algebroid consists in a differential module structure on $\Gamma(T^*(M)).$
Indeed, one has
\bea
\label{mod}
\{\a,f\b\}=f\{\a,\b\}+(\lc_{\po\a}f)\b.
\eea
{\lemma \label{hom1}
\bea
\label{hom}
[\po \a,\po \b]=\po \{\a,\b\}
\eea
}
\\
{\bf Proof~} For $\a,\b$ exact the statement follows from the standard fact in
Poisson geometry that the commutator of Hamiltonian vector fields is a
Hamiltonian vector field corresponding to the Poisson bracket of
Hamiltonians
$$[X_f,X_g]=X_{\{f,g\}}$$
due to the fact that $\po df =X_f.$ Then, locally each $1$-form can be
represented as a sum of expressions of the type  $fdg,$
and it is sufficient to prove the claim for $\a=f_1df_2,~\b=g_1dg_2$ due to
the local nature of this relation.
One has
\bea
[\po f_1df_2,\po g_1dg_2]&=&f_1g_1[\po df_2,\po
dg_2]+f_1(\lc_{\po\a}g_1)\po\b-g_1(\lc_{\po\b}f_1)\po\a\nn\\
\po\{f_1df_2,g_1dg_2\}&=&\po(f_1g_1\{df_2,dg_2\}+
f_1(\lc_{\po\a}g_1)\b-g_1(\lc_{\po\b}f_1)\a)\nn\\
&=&f_1g_1\po\{df_2,dg_2\}+
f_1(\lc_{\po\a}g_1)\po\b-g_1(\lc_{\po\b}f_1)\po\a\nn
\eea
Then returning to the formula (\ref{hom}) for exact $1$-forms one obtains the general
statement $\square$
{\lemma The bracket (\ref{p}) satisfies the
Jacobi identity. }
\\
{\bf Proof~} We proceed as in the previous lemma starting from consideration
of exact $1$-forms.
The Jacobi identity for forms $\a=df,~\b=dg,~\g=dh$
is equivalent to the Jacobi identity for the bracket of functions:
\bea
J(df,dg,dh)&=&\{\{df,dg\}dh\}+\{\{dg,dh\}df\}+\{\{dh,df\}dg\}=d
J(f,g,h)\nn
\eea
where on the right-hand side the Jacobi expression is
taken with respect to the Poisson bracket on functions. To finish the proof
one needs to show that the Jacobi expression $J(\a,\b,\g)$ is
tri-linear subject to multiplication by functions. Indeed,
\bea
J(\a,\b,f\g)&=&\{\{\a,\b\}f\g\}+\{\{\b,f\g\}\a\}+\{\{f\g,\a\}\b\}\nn\\
&=&f\{\a,\b\},\g\}+(\lc_{\po \{\a,\b\}}f)\g+f\{\{\b,\g\},\a\}-(\lc_{\po
\a}f)\{\b,\g\}\nn\\
&+&(\lc_{\po \b}f)\{\g,\a\}-(\lc_{\po \a}(\lc_{\po \b}f))\g+f\{\{\g,\a\},\b\}-
(\lc_{\po \b}f)\{\g,\a\}\nn\\
&-&(\lc_{\po \a}f)\{\b,\g\}+(\lc_{\po \b}(\lc_{\po \a}f))\g\nn\\
&=&fJ(\a,\b,\g)\nn
\eea
where we used the property (\ref{mod}) and the result of lemma \ref{hom1}
$\square$

\paragraph{\fbox{Extension to $n$-forms}}
In fact the Lie bracket (\ref{p}) can be extended to the graded Lie bracket
on the de Rham complex. To do this we first observe, that it can be extended to the
exterior algebra of $1$-forms $\bigwedge^*\Gamma\Omega^1(M)$ as follows.
 Let us
introduce the grading $deg \o=d =|\o|-1$ for
$\o \in \bigwedge^{d+1}\Gamma\Omega^1(M).$ Further we
denote by $d_i$ the degree $deg\o_i.$
Now one can use the following formula to define the sought extension of (\ref{p})
\bea
\label{ext}
\{\o_1,\o_2\wedge\o_3\}=\{\o_1,\o_2\}\wedge\o_3+
(-1)^{d_1(d_2+1)}\o_2\wedge\{\o_1,\o_3\}
\eea
Imposing the graded antisymmetry
$$\{\o_1,\o_2\}=-(-1)^{d_1d_2}\{\o_2,\o_1\}$$
one obtains
\bea
\{\o_1\wedge\o_2,\o_3\}=(-1)^{d_3(d_2+1)}\{\o_1,\o_3\}\wedge\o_2+
\o_1\wedge\{\o_2,\o_3\}\nn
\eea
This bracket satisfies the graded Jacobi identity:
\bea
0&=&J_{gr}(\o_1,\o_2,\o_3)\nn\\
&=&(-1)^{d_1 d_3}\{\{\o_1,\o_2\},\o_3\}
+(-1)^{d_2 d_1}\{\{\o_2,\o_3\},\o_1\}
+(-1)^{d_3 d_2}\{\{\o_3,\o_1\},\o_2\}.\nn
\eea
We can proceed by induction on the degree of $\o_3$ due to the symmetry of
the Jacobi relation. In fact,
\bea
J(\o_1,\o_2,\o_3\wedge\a)=(-1)^{d_1}J(\o_1,\o_2,\o_3)\wedge\a+
(-1)^{d_2(d_3+1)}\o_3\wedge J(\o_1,\o_2,\a)\nn
\eea
and hence it reduces to the Jacobi identity for lower order elements.

In fact the bracket on the exterior algebra $\bigwedge^*\Gamma\Omega^1(M)$
descends correctly to the de Rham complex with the same grading. To verify
this it is sufficient to prove the next formula:
\bea
\label{leib}
\{\o_1,f\o_2\wedge\o_3\}=\{\o_1,\o_2\wedge f\o_3\}.
\eea
We shall restrict ourselves to the decomposable case $\o_1=\a_1\w\ldots\w\a_k$ and
$\o_2=\b_1\w\ldots\w\b_l$ and demonstrate that the
the following equation holds for all $i$ and $j$:
\bea
\label{leib1}
\{\a_1\w\ldots\w\a_k,\b_1\w\ldots\w f\b_i\w\ldots\w\b_l\}=
\{\a_1\w\ldots\w\a_k,\b_1\w\ldots\w f\b_j\w\ldots\w\b_l\}
\eea
Consider the left-hand side of \eqref{leib1}
\bea
&&\sum_{m}(-1)^{m+i}\{\a_m,f\b_i\}\w\a_1\w\ldots\w\widehat{\a_m}\w\ldots\w\a_k\w
\b_1\w\ldots\w\widehat{\b_i}\w\ldots\w\b_l\nn\\
&+&f\sum_{m,k\neq i}(-1)^{m+k}\{\a_m,\b_k\}\w\a_1\w\ldots\w\widehat{\a_m}\w\ldots\w\a_k\w
\b_1\w\ldots\w\widehat{\b_k}\w\ldots\w\b_l\nn\\
&=&\sum_{m}(-1)^{m+i}\lc_{\po\a_m}f\b_i\w\a_1\w\ldots\w\widehat{\a_m}\w\ldots\w\a_k\w
\b_1\w\ldots\w\widehat{\b_i}\w\ldots\w\b_l\nn\\
&+&f \sum_{m}(-1)^{m+i}\{\a_m,\b_i\}\w\a_1\w\ldots\w\widehat{\a_m}\w\ldots\w\a_k\w
\b_1\w\ldots\w\widehat{\b_i}\w\ldots\w\b_l\nn\\
&+&f\sum_{m,k\neq i}(-1)^{m+k}\{\a_m,\b_k\}\w\a_1\w\ldots\w\widehat{\a_m}\w\ldots\w\a_k\w
\b_1\w\ldots\w\widehat{\b_k}\w\ldots\w\b_l\nn\\
&=&f\{\o_1,\o_2\}+\lc_{\po\o_1}f\w\o_2;
\eea
where
\bea
\label{lie2}
\lc_{\po\o_1}f=\sum_{m}(-1)^{m}(\lc_{\po\a_m}f)\a_1\w\ldots\w\widehat{\a_m}
\w\ldots\w\a_k.
\eea
Now it is enough to notice that the resulting formula does not
depend on $i.$

\paragraph{\fbox{Reduction}}
The next crucial property of the considered construction that makes $\Omega^*(M)$ a DGLA is
relation between the graded Lie algebra structure on $\Omega^*(M)$ and the
de Rham differential, exactly it is that $d$
differentiates the Lie bracket introduced above.
{\lemma
\bea
\label{dif}
d\{\o_1,\o_2\}=\{d\o_1,\o_2\}+(-1)^{d_1}\{\o_1,d\o_2\}
\eea
}
\\
{\bf Proof~} We proceed by induction on the degree of $\o_2.$
We take a decomposable element $\o_2\w\o_3$ and suppose that
\bea
d\{\o_1,\o_2\}&=&\{d\o_1,\o_2\}+(-1)^{d_1}\{\o_1,d\o_2\};\nn\\
d\{\o_1,\o_3\}&=&\{d\o_1,\o_2\}+(-1)^{d_1}\{\o_1,d\o_3\}.\nn
\eea
Then one needs to prove that
\bea
d\{\o_1,\o_2\w\o_3\}=\{d\o_1,\o_2\w\o_3\}+
(-1)^{d_1}\{\o_1,d(\o_2\w\o_3)\}.\nn
\eea
On the left-hand side we have:
\bea
&&d(\{\o_1,\o_2\}\w\o_3+(-1)^{d_1(d_2+1)}\o_2\w\{\o_1,\o_3\})\nn\\
&=&\{d\o_1,\o_2\}\w\o_3+(-1)^{d_1}\{\o_1,d\o_2\}\w\o_3+
(-1)^{d_1+d_2+1}\{\o_1,\o_2\}\w d\o_3\nn\\
&+& (-1)^{d_1(d_2+1)}d\o_2\w\{\o_1,\o_3\}+(-1)^{(d_1+1)(d_2+1)}
(\o_2\w\{d\o_1,\o_3\}+(-1)^{d_1}\o_2\w\{\o_1,d\o_3\}).\nn
\eea
On the right-hand side we have:
\bea
&& \{d\o_1,\o_2\}\w\o_3+
(-1)^{(d_1+1)(d_2+1)}\o_2\w\{d\o_1,\o_3\}\nn\\
&+& (-1)^{d_1}(\{\o_1,d\o_2\}\w\o_3+
(-1)^{d_1(d_2+2)}d\o_2\w\{\o_1,\o_3\}\nn\\
&+& (-1)^{d_2+1}\{\o_1,\o_2\}\w
d\o_3+(-1)^{d_2+1+d_1(d_2+1)}\o_2\w\{\o_1,d\o_3\}).\nn
\eea
Comparing both sides one obtains the lemma $\square$
{\lemma The Lie algebra structure on
$\Omega^*(M)$ can be pushed down to its de Rham cohomology}
\\
{\bf Proof~~} Let $\o_1,\o_2$ be closed then $\{\o_1,\o_2\}$ is also closed;
let $\o_1$ be closed and $\o_2$ - exact, then $\{\o_1,\o_2\}$ is exact due
to (\ref{dif}) $\square$

{\lemma The induced Lie algebra structure on $H^*(M)$ of a Poisson manifold
is trivial.}
\\
{\bf Proof~~} It is obvious for $1$-forms. Indeed, let $\a,\b$ be closed:
$d\a=0,d\b=0.$ Then
$$\lc_{\po\a}\b=d\imath_{\po\a}\b=d\pi(\a,\b)$$ and
\begin{equation}\label{eq5 1/4}\{\a,\b\}=d\pi(\a,\b)\end{equation}
so, it is exact.

Let us next consider closed forms $\a,\o$ of degrees $1$ and $k$
respectively. There is a simplification of the formula (\ref{ext}) for this
case:
\bea
\{\a,\o\}=\lc_{\po \a}\o=d(\imath_{\po\a}\o).\nn
\eea
Indeed, on both sides one has a differentiation and the
equality fulfills when $\o$ is a
$1$-form.

Now observe that the right hand side of the last formula can be written down in
the form
\begin{align}
\label{eq5 1/2}
\{\o_1,\,\o_2\}&=d\widetilde\pi(\o_1,\,\o_2),\\
\intertext{where $\widetilde\pi$ is a bilinear operation, sending two forms
$\o_1,\,|\o_1|=k$ and $\o_2,\,|\o_2|=l$ to
a $(k+l-2)$-form}
\label{eq5 3/4}
\widetilde\pi(\o_1,\o_2)=\sum_{p,q}\sum_{1\le i\le k,\,1\le
j\le
l}(-1)^{k+j-1}&\pi(\alpha_i,\,\beta_j)
\wedge\alpha^p_1\dots\wedge\widehat\alpha^p_i\wedge\dots\wedge\widehat\beta^q_j
\wedge\dots\wedge\beta^q_l,\\
\intertext{where, as usually,\ \ $\widehat{}$\ \ over an element means that this element is omitted and}
\notag
\o_1=\sum_p\alpha^p_1\wedge&\dots\wedge\alpha^p_k,\\
\notag
\o_2=\sum_q\beta^q_1\wedge&\dots\wedge\beta^q_l.
\end{align}
As before, $\alpha^p_i,\ \beta^q_j$ are 1-forms. Observe, that
$$\widetilde\pi(\o_1,\o_2)=\pi\vdash(\o_1\w\o_2)-
(\pi\vdash\o_1)\w\o_2-(-1)^{|\o_1|}\o_1\w(\pi\vdash\o_2),$$
where $\pi\vdash$ denotes the internal multiplication of a differential form
by the bivector $\pi$. In particular, it follows that $\widetilde\pi$ is a
bilinear with respect to the multiplication by functions and
differentiation with respect to either argument (which can also be proved
independently directly from the formula \eqref{eq5 3/4}).

We shall prove, that the formula \eqref{eq5 1/2} holds in the general case,
i.e that {\em $\o_1,\ \o_2$ are
closed forms, then their Poisson bracket is an exact form, defined by formulas \eqref{eq5 1/2}, \eqref{eq5 3/4}\/}.

To this end,
let us observe that the right hand side of the formula \eqref{eq5 3/4} (and
consequently that of \eqref{eq5 1/2}) is well-defined, that is it doesn't depend on the way
one decomposes the forms $\alpha,\ \beta$ into the wedge-product of 1-forms.
Second, it is enough to prove formula \eqref{eq5 1/2} locally, i.e. on an
arbitrary open domain in the manifold. Indeed, the expressions on both its sides are
well-defined everywhere on the manifold and don't depend
on the choice of the decompositions of $\alpha,\ \beta$ into $1$-forms. Finally, recall, that locally all closed forms are
exact, and hence they are representable as the sum of products of
closed (and even exact) forms (one has $\o=d\o'=d(\sum_{I=(i_1,\dots,i_n)}f_Idx^{i_1}\w\dots\w dx^{i_n})=\sum_Idf_Idx^{i_1}\w\dots\w dx^{i_n}$). Now the conclusion follows from the formula
\eqref{eq5 1/4}.$\square$
{\Rem\rm
In effect, one can prove the following formula (pointed out by I.Nikonov):
\begin{equation}
\label{eq5 4/5}
\{\o_1,\o_2\}=d\widetilde\pi(\o_1,\,\o_2)-
\widetilde\pi(d\o_1,\,\o_2)-(-1)^{|\o_1|}\widetilde\pi(\o_1,\,d\o_2).
\end{equation}
Here $\o_1,\ \o_2$ are arbitrary de Rham forms on the manifold and
$\widetilde\pi$ is the operation, defined above by the formula \eqref{eq5
3/4}.

In order to prove \eqref{eq5 4/5} just observe, that its right hand side is a skew-symmetric bilinear form of degree $-1$ on the de
Rham complex of the Poisson manifold, verifying the Leibnitz rule with respect to the
multiplication of the forms (i.e. the formula \eqref{leib}
holds for $\{\,,\,\}$ replaced with the expression from \eqref{eq5 4/5}),
and differentiated by the de Rham differential (i.e. verifying the formula
\eqref{dif}). Finally, observe that it coincides with the bracket $\{\,,\,\}$ when both forms have
degree 1.}

\section{$L_\infty$-\/algebras and $L_\infty$\/-morphisms}
In this section we recall the basic definitions and constructions of the
homotopy Lie algebras ($L_\infty$-algebras) and $L_\infty$ morphisms, which are crucial for our
treatment of the formality (see section 1). First, let us recall the definition of
$L_\infty$-algebras. In what follows, we use the notation and signs from the
paper \cite{PenWel}.

\paragraph{\fbox{The free cocommutative coalgebra of a graded space}}
Let $L$ be a graded vector space. Consider the external algebra, associated to
$L$:
\begin{equation}
\label{eq5}
\wedge^*L=\bigoplus_{i\ge1}L^{\otimes i}\Bigl/\{a\otimes b+(-1)^{|a||b|}b\otimes
a\}.
\end{equation}
Here and below $|a|$ will denote the degree of an element $a$ in $L$.
Recall, that one can introduce double grading on $\Lambda^*L$,
putting $\textrm{bideg}(a_1\wedge a_2\wedge\dots\wedge a_n)=(\sum |a_i|,n),\
a_i\in L$. For an element of bi-degree $(m,n)$ we define its total
degree to be equal to the difference $n-m$. We shall denote the second degree of an element $\alpha$ by
$|\alpha|'$ (in contrast with its first degree, denoted by $|\a|$).

One uses
the following formula as the definition of comultiplication $\nabla$ in $\wedge^*L$
($a_i\in L$ are arbitrary elements):
\begin{equation}
\label{eq6}
\nabla(a_1\wedge\dots\wedge a_n)=\sum_{i=1}^{n-1}\sum_{\sigma\in
Sh_{k,l}}(-1)^{\sigma+kl}a_{\sigma(1)}\wedge\dots\wedge a_{\sigma(k)}\otimes
a_{\sigma(k+1)}\wedge\dots\wedge a_{\sigma(n)}.
\end{equation}
Here $Sh_{k,l}$ stands for the subset of all $(k,l)$-shuffles in the group
of permutations of $k+l$ elements, that is of the permutations, preserving
the order of the first $k$ and the last $l$ elements. And the sign $(-1)^\sigma$
of the shuffling $\sigma$, is determined by the following
rule: every time we change the order of two elements $a_i$
and $a_j$ in the sequence $a_1,\dots,a_n$, we multiply it by $|a_i||a_j|+1$.
Equipped with this comultiplication $\wedge^*L$ is the free cocommutative
coalgebra, (without counit) (co)generated by $L$.

Clearly, the comultiplication $\nabla$ respects the double grading and hence the total grading in $\wedge^*L$, so
that the latter space becomes a (bi)graded cocommutative coalgebra.

Now one can introduce bi-degrees for all the maps from
$\Lambda^*L$ to itself, as (minus) the difference of degrees of a
(bi-homogeneous) element and its image. Of course, this is not always
correctly defined (i.e. it may depend on the element). We would
say that the map is (bi-)homogeneous, if it has bi-degree. As in
the case of elements, we denote the second degree of a map
$f$ by $|f|'$, reserving the symbol
$|f|$ for the first (usual) degree of the map $f$.

As before, one defines the total degree of a map as the difference of its
first and second degrees. Observe, that a map, which is homogeneous with respect to
this grading, is not necessarily bihomogeneous (while the contrary is true).

\paragraph{\fbox{Graded coderivatives}}
One can give the following important definition
\begin{df}
Let $(C,\ \nabla_C)$ be a bigraded coalgebra with comultiplication $\nabla_C$.
A bi-homogeneous map $D:C\to C$ of bidegree $(|D|,\,|D|')$
is called {\em (bi-homogeneous) coderivative\/}, if the following diagram
is commutative up to a sign, depending on the bidegrees of elements:
\begin{equation}
\label{eq7}
\begin{CD}
C @>{\nabla_C}>>{C\otimes C}\\
@VVDV            @VV{D\otimes1+1\otimes D}V\\
C @>{\nabla_C}>>{C\otimes C.}
\end{CD}
\end{equation}
Or, in more explicit terms,
$$
\nabla(D\alpha)=D\alpha_{(1)}\otimes\alpha_{(2)}+(-1)^{|D||\alpha_{(1)}|+|D|'|\alpha_{(1)}|'}\alpha_{(1)}\otimes D\alpha_{(2)},
$$
for all $\alpha$ in $C$. Here we use the Sweedler's notation for
comultiplication: $\nabla(\alpha)=\alpha_{(1)}\otimes\alpha_{(2)}$.
\end{df}
\begin{prop}
\label{prop1}
Let $D_1$ and $D_2$ be two bi-homogeneous coderivatives of
$\Lambda^*V$. Then their bi-graded commutator
\begin{equation}
\label{commut}
[D_1,\,D_2]=D_1D_2-(-1)^{|D_1||D_2|+|D_1|'|D_2|'}D_2D_1
\end{equation}
is also a bi-homogeneous
coderivative. Its bi-degree is equal to the (element-wise) sum of
the bi-degrees of $D_1$ and $D_2$.
\end{prop}
\begin{proof}
Evident.
\end{proof}

By the virtue of freeness of the coalgebra $\wedge^*L$, any coderivative on it is uniquely
determined by its ``values on cogenerators", that is by a collection of linear
mappings $l_n:\wedge^nL\to L$ (see \cite{Swe}). If the total degree of $l$ is equal to $p$, we conclude, that
the $n$-th map
in this set, $l_n$, sends elements $a_1\wedge\dots\wedge a_n\in\wedge^nL$ to the
elements of degree $|a_1|+\dots+|a_n|-n+p+1$ in $L$, due to the introduced
grading in $\wedge^*L$. And the converse is also true: if $\{f_k:\Lambda^kL\to\Lambda^1L\}$ is an arbitrary collection of linear
maps, such that the first degree of $f_k$ is equal to $m_k$, then one can in a unique way extend $\{f_k\}$ to a coderivative $F:\Lambda^*V\to\Lambda^*V$, equal to the sum
of bi-homogeneous coderivatives of bi-degrees $(m_k,k-1)$. (Here the phrase ``$F$ extends $f$" means
that $F_{|\Lambda^iV}=f_i$.) This is an immediate consequence of freeness of $\Lambda^*V$ (see Sweedler's
book). More explicitly, one uses following formula as the definition of $F$:
\begin{equation}\label{eqF}
F(a_1\wedge\dots\wedge
a_n)=\sum_{k<n}\sum_{\sigma\in\Sigma_n}(-1)^{\sigma}\frac{1}{k!(n-k)!}f_k(a_{\sigma(1)}\wedge\dots\wedge
a_{\sigma(k)})\wedge a_{\sigma(k+1)}\wedge\dots\wedge
a_{\sigma(n)}
\end{equation}
Here the sign $(-1)^{\sigma}$ is defined as above.

One can check, that the commutator of two bi-homogeneous
coderivatives $\Phi$ and $\Psi$ of bidegrees $(|\varphi|,\,k-1)$ and
$(|\psi|,\,l-1)$, extending the maps $\varphi:\Lambda^kL\to\Lambda^1L$ and
$\psi:\Lambda^lL\to\Lambda^1L$ coincides with the extension of the map
$[\varphi,\,\psi]:\Lambda^{k+l-1}L\to\Lambda^1L$, given by the formula
\begin{equation}
\label{eqnew}
\begin{split}
&[\varphi,\,\psi](a_1\wedge\dots\wedge a_{k+l-1})=\\
                &\quad\sum_{\sigma\in\Sigma_{k+l-1}}\frac{1}{l!(k-1)!}\varphi(\psi(a_{\sigma(1)}\wedge\dots\wedge a_{\sigma(l)})\wedge a_{\sigma(l+1)}\wedge\dots\wedge
a_{\sigma(k+l-1)})-\\
                &\quad(-1)^\epsilon\sum_{\sigma\in\Sigma_{k+l-1}}\frac{1}{(l-1)!k!}\psi(\varphi(a_{\sigma(1)}\wedge\dots\wedge a_{\sigma(k)})\wedge a_{\sigma(k+1)}\wedge\dots\wedge
a_{\sigma(k+l-1)}),
\end{split}
\end{equation}
where $\epsilon={|\varphi||\psi|+(k-1)(l-1)}$.

The space of (bi-graded) coderivatives of $\wedge^*L$, endowed with the commutator
$[,\,]$, turns into a {\em bi-\/}graded Lie algebra. However, one can turn
it into the graded Lie algebra with respect to the total degree simply by
changing slightly the sign in \eqref{commut}, namely
\begin{equation}
\label{defcom}
\{D_1,\,D_2\}=(-1)^{|D_1||D_2|'}[D_1,\,D_2].
\end{equation}

An alternative way of looking at this formula is as follows. Consider a
``deformation" $\bar D$ of a bihomogeneous map $D:\wedge^*L\to\wedge^*L$,
defined by the formula $\bar D(\a)=(-1)^{|D||\a|'}D(\a)$. Then
$\overline{D_1D_2}=(-1)^{|D_1||D_2|'}\bar{D_1}\bar{D_2}$. Indeed,
\begin{equation*}
\begin{split}
\bar{D_1}(\bar{D_2}(\a))&=(-1)^{|D_2||\a|'}\bar{D_1}(D_2(\a))\\
                      &=(-1)^{|D_1|(|D_2|'+|\a|')+|D_2|'|\a|}D_1(D_2(\a))=(-1)^{|D_1||D_2|'}\overline{D_1D_2}(\a).
\end{split}
\end{equation*}
Observe, that if $D$ is a coderivative, then $\bar D$ is {\em not\/} a bi-graded coderivative
for the given diagonal $\nabla$ (the signs rule is violated). In effect, it is a {\em graded\/} coderivative (i.e. a coderivative with respect to the total degree), for a ``deformed"
comultiplication $\bar\nabla$, given by the formula
$\bar\nabla(\a)=(-1)^{|\a_{(1)}||\a_{(2)}|'}\a_{(1)}\otimes\a_{(2)}$, where
$\nabla(\a)=\a_{(1)}\otimes\a_{(2)}$. This is verified by a direct
computation:
\begin{equation*}
\begin{split}
\bar\nabla(\bar D(\a))&=(-1)^{|D||\a|'}\bar\nabla(D(\a))=\\
                    &=(-1)^{|D||\a|'}\Bigl((-1)^{(|D|+|\a_{(1)}|)|\a_{(2)}|'}D(\a_{(1)})\otimes\a_{(2)}\\
                    &\quad+(-1)^{|D||\a_{(1)}|+|D|'|\a_{(1)}|'+|\a_{(1)}|(|D|'+|\a_{(2)}|')}\a_{(1)}\otimes
D(\a_{(2)})\Bigr)\\
                    &=(-1)^{|D||\a|'}\Bigl((-1)^{|D|(|\a_{(1)}|'+|\a_{(2)}|')+|\a_{(1)}||\a_{(2)}|'}\bar D(\a_{(1)})\otimes\a_{(2)}\\
                    &\quad+(-1)^{|D||\a_{(1)}|+|D|'|\a_{(1)}|'+|\a_{(1)}|(|D|'+|\a_{(2)}|')+|D||\a_{(2)}|'}\a_{(1)}\otimes\bar
D(\a_{(2)})\Bigr)\\
                    &=(-1)^{|\a_{(1)}||\a_{(2)}|'}(\bar
D(\a_{(1)})\otimes\a_{(2)}+(-1)^{(|D|+|D|')(|\a_{(1)}|+|\a_{(1)}|')}\a_{(1)}\otimes\bar
D(\a_{(2)}))\\
                    &\qquad\qquad\qquad\qquad=(\bar D\otimes 1+(-1)^{(|D|+|D|')(|\a_{(1)}|+|\a_{(1)}|')}1\otimes\bar
D)\bar\nabla(\a).
\end{split}
\end{equation*}
Now, if $[A,\,B]_{tot}$ denotes the commutator of two maps
$A,B:\wedge^*L\to\wedge^*L$ with respect to their {\em total\/} degree, then
one has
\begin{equation}
\label{eqsigns}
\begin{split}
[\bar{D_1},\,\bar{D_2}&]_{tot}=\bar{D_1}\bar{D_2}+(-1)^{(|D_1|+|D_1|')(|D_2|+|D_2|')}\bar{D_2}\bar{D_1}\\
                           &=(-1)^{|D_1||D_2|'}\overline{D_1D_2}+(-1)^{|D_2||D_1|'+(|D_1|+|D_1|')(|D_2|+|D_2|')}\overline{D_2D_1}\\
                           &=(-1)^{|D_1||D_2|'}(\overline{D_1D_2}+(-1)^{|D_1||D_2|'+|D_2||D_1|'+(|D_1|+|D_1|')(|D_2|+|D_2|')}\overline{D_2D_1})\\
                           &=(-1)^{|D_1||D_2|'}(\overline{D_1D_2}+(-1)^{|D_1||D_2|+|D_1|'|D_2|'}\overline{D_2D_1})=\overline{\{D_1,\,D_2\}}
\end{split}
\end{equation}

\paragraph{\fbox{$L_\infty$\/-algebras and morphisms}}
The definition of an $L_\infty$-algebras can now be formulated in a very
concise way
\begin{df}
\label{defL}
One says, that the graded space $L$ is an $L_\infty$-\/algebra,
if its free coalgebra $\wedge^*L$ is equipped with a degree-1 (i.e. its total degree should be equal to 1) bigraded coderivative $l$
verifying the relation $\{l,\,l\}=0$. Equivalently, in the view of equation \eqref{eqsigns},
one can say, that the graded coderivative $\bar l$ verifies the equation $(\bar{l})^2=0$.
\end{df}

The equations, appearing in this definition, can be interpreted as an
infinite collection of quadratic equations on the components $l_i$ of $l$.
Here are the first three equations of this collection:
\begin{gather*}
l_1(l_1(a))=0;\\
l_1(l_2(a\wedge b))=l_2(l_1(a)\wedge b)-(-1)^{|a|}l_2(a\wedge l_1(b));\\
\begin{split}
l_1(l_3(a\wedge b\wedge c))&+l_3(l_1(a)\wedge b\wedge c)+(-1)^{|a|}l_3(a\wedge l_1(b)\wedge c)\\
                           &\qquad\qquad+(-1)^{|a|+|b|}l_3(a\wedge b\wedge l_1(c))\\
                           &=\frac12\Bigl((-1)^{|a||c|}l_2(l_2(a\wedge b)\wedge c)+(-1)^{|a|(|b|+|c|)}l_2(l_2(b\wedge c)\wedge a)\\
                           &\qquad\qquad+(-1)^{(|a|+|b|)|c|}l_2(l_2(c\wedge a)\wedge
b)\Bigr).
\end{split}
\end{gather*}
The first and the second of these equalities mean, that $l_1$ is a
square-zero differential and $l_2$ -- an anti-symmetric bracket on the space
$L$, so that $l_1$ is a differentiation of $l_2$. The third equation
implies, that the Jacobi identity on $L$ holds up to a homotopy, moreover,
the chain homotopy connecting it (the Jacobi formula) to zero is $l_3$.

On the whole, the equations, verified by the mappings $l_i$ can be written
down in the following way:
\begin{equation}
\sum_{\substack{k+l=n+1\\ \sigma\in Sh_{k,l-1}}}(-1)^{\sigma+(k-1)l}l_l(l_k(a_{\sigma(1)}\wedge\dots\wedge a_{\sigma(k)})\wedge\dots\wedge
a_{\sigma(n)})=0.
\end{equation}
Or, in a more conceptual form,
\begin{equation*}
d(l_n)=\sum_{i=1}^{[\frac{n}{2}]}\frac12l_i\circ(l_{n-i}\times1).
\end{equation*}
Here we have denoted the differential $l_1$ in $L$ (see above) by a more traditional letter $d$, $[x]$ denotes the integer part of the number $x$, $d(l_n)=d\circ l_n+(-1)^{n-2}l_n\circ d$,
and $l_i\times 1\eqdef(l_i\otimes 1)\circ\nabla$. In brief, $l_n$ is a chain homotopy connecting certain linear
combination of $l_i,\ i<n$ with zero.

In particular, one sees from these formulas, that {\em every differential graded
Lie algebra is $L_\infty$\/-algebra}. It is enough to put
\begin{align*}
l_1&=d \ \ \quad{\text{-- differential in $L$,}}\\
l_2&=[\,,] \quad{\text{-- the Lie bracket in $L$, and}}\\
l_i&=0,\qquad i\ge2.
\end{align*}

One can define two different types of natural maps between $L_\infty$\/-algebras.
First is the so-called {\em strict} $L_\infty$\/-maps, i.e. the maps,
commuting with all the mappings $l_n$. The second, and the most important
type is the {\em strong homotopy} morphisms or $L_\infty$\/-morphisms. By
definition, an $L_\infty$\/-morphism from an $L_\infty$\/-algebra
$L$ to an $L_\infty$\/-algebra $L'$ is a (total) degree-0 linear map $F:\wedge^*L\to\wedge^*L'$,
commuting with comultiplications and (co)differentials $l,\ l'$ in $\wedge^*L$ and $\wedge^*L'$.

As before, due to the freeness of $\wedge^*L'$, any homomorphism to this
coalgebra from a coalgebra $C$ is uniquely determined by its ``zero stage",
i.e. by the linear map $C\to L$, which it defines. In our situation, this
reduces to a collection of maps  $f_n:\wedge^nL\to L',\ n\ge1$, each changing the
(first) degree by $-(n-1)$. The condition that the total map, made up from
all the $f_i$ commutes with the (co)differentials, takes the form of an infinite
series of equations on the maps $f_i$. The first two equations of this
series are as follows (c.f. \eqref{eq2},\eqref{eq3}):
\begin{align}
\label{eq8}
l'_1(f_1(a))&=f_1(l_1(a)),\\
\label{eq9}
l'_1(f_2(a\wedge b))-f_2(l_1(a)\wedge b)-(-1)^{|a|}&f_2(a\wedge l_1(b))=f_1(l_2(a\wedge b))-l'_2(f_1(a)\wedge
f_1(b)).
\end{align}

The equation \eqref{eq8} shows that the map $f_1$ should commute with the
differentials $l_1,\ l'_1$ in $L,\ L'$, and the equation \eqref{eq9} means
that $f_2$ is the chain homotopy, which makes the map $f_1$ a homomorphism
of Lie algebras on the corresponding cohomology. The general equation from
the definition of $L_\infty$\/-morphisms can be written down in the following
form:
\bea
\notag
df_n&=&\sum_{k=2}^n\frac{1}{k!}\ \sum_{i_1+\dots+i_k=n}\pm
l_k\circ(f_{i_1}\times\dots\times f_{i_k})\\
\notag
&+&\sum_{k=1}^{n-1}\pm\frac12f_k\circ(l_{n-k+1}\times1).
\eea
Here, as above, we denote by $f_1\times f_2$ the composition $(f_1\otimes f_2)\circ\nabla$
(the ``external product" of more than $2$ terms is defined by induction).

The signs in this formula depend on the dimensions of the maps, and on their
order. Since below we shall consider only the case when both
$L_\infty$\/-algebras are in effect (graded, differential) Lie algebras, let
us give the precise formulas for the equations, verified by an
$L_\infty$\/-morphism just in this case:
\begin{equation}
\label{eq10}
\begin{split}
&df_n(a_1\wedge\dots\wedge a_n)=\\
&=\sum_{i=1}^n(-1)^{n-1+\varepsilon_i}f(a_1\wedge\dots\wedge da_i\wedge\dots\wedge
a_n)\\
&+\frac{1}{2}\sum_{i+j=k}\sum_{\sigma\in
Sh_{i,j}}(-1)^{\sigma+(i-1)\varepsilon_i}\{f_i(a_{\sigma(1)}\wedge\dots\wedge a_{\sigma(i)}),\,f_j(a_{\sigma(i+1)}\wedge\dots\wedge
a_{\sigma(n)})\}\\
&+\frac12\sum_{1\le i<j\le
n}(-1)^{(i-1)\varepsilon_{i-1}+(j-1)\varepsilon^i_{j-1}}f_{n-1}(\{a_i,\,a_j\}\wedge a_1\wedge\dots\wedge\widehat a_i\wedge\dots\wedge\widehat a_j\wedge\dots\wedge
a_n).
\end{split}
\end{equation}
As usually, the ``hat" over an entry means that the entry is omitted in the
formula, $Sh_{i,j}$ denotes the collection of all $(i,\,j)$\/-shuffles,
the sign $(-1)^\sigma$ is defined as above and
$$\varepsilon_i=-i+\sum_{k=1}^i|a_k|,\ \varepsilon^j_i=\varepsilon_i-|a_j|+1.$$

\section{Main theorem}
The purpose of this section is to prove the formality of the DG Lie algebra,
associated to a Poisson manifold. To this end we shall construct an
$L_\infty$\/-quasi-isomorphism (i.e. an $L_\infty$\/-map with quasi-isomorphic first stage) between this Lie algebra and its cohomology
Lie algebra. More accurately, in order to speak about the Lie algebras in this case, one should
first change grading in $\Omega^*(\mc)$. In fact, the bracket, introduced in
\S2 sends a couple of elements of degrees $k$ and $l$ in $\Omega^*(\mc)$ to
an element of degree $k+l-1$ and not $k+l$. So, we pass to the suspension $\Sigma\Omega^*(\mc)$ of $\Omega^*(\mc)$,
Here for a graded space $V$, $\Sigma V$ is the graded space, defined by
$$
(\Sigma V)_i=V_{i+1}.
$$
(In particular, $\Sigma\Omega^*(\mc)_{-1}=C^\infty(\mc)$.) It is easy to
check, that the new grading is respected by the Cartan bracket on
$\Omega^*(\mc)$.
Similarly, one should replace $H(\mc)$ by its suspension $\Sigma H(\mc)$.

In order to prove the formality of the Lie algebra $\Sigma\Omega^*(\mc)$ for
a Poisson manifold $\mc$, we shall need some special properties of the free
coalgebra, generated by $\Sigma\Omega(\mc)$. These properties result from the
usual constructions, related to the differential forms on a smooth manifold, such as the Cartan calculus and external multiplication of forms. We have collected these results in a
separate subsection.

\subsection{Coderivatives in $\wedge^*\Sigma\Omega(\mc)$}
Let $L=\Sigma\Omega(\mc)$ -- the suspension of the algebra of de Rham forms for
some manifold $\mc$ (see above). We shall denote by $s$ the evident degree $-1$
isomorphism $\Omega(\mc)\stackrel{s}{\to}\Sigma\Omega(\mc)$, in particular
$s\alpha$ will denote the element of $\Sigma\Omega(\mc)$, corresponding to a form
$\alpha\in\Omega(\mc)$.

Below we give few geometric
examples of coderivatives of $\Lambda^*L$ in this case.
\paragraph{\fbox{Cartan algebra}}
Let $X$ be a vector field. We can associate to it two differentiations of $\Omega(\mc)$,
Lie derivative and the internal multiplication by $X$. We shall denote them by $\lc_X$ and $\imath_X$ respectively.
The degree of $\lc_X$ is equal to $0$ and of $\imath_X$ --- to $-1$.

We shall associate to $\lc_X$ and $\imath_X$ the
maps $\tilde\lc_x,\ \tilde\imath_X:\Sigma\Omega(\mc)\to\Sigma\Omega(\mc)$
as follows:
$$\tilde\lc_X(s\alpha)=s\lc_X(\alpha),\quad\tim_X(s\alpha)=-s\imath_X(\alpha).$$
The degrees of these maps are also equal to $0$ ad $-1$.

These two maps can be extended to coderivatives of
$\Lambda^*\Sigma\Omega(\mc)$ of bi-degrees $(0,0)$ and $(0,-1)$ respectively (see the remark, following the proposition \ref{prop1}). We shall denote these
coderivatives by $L_X,\ I_X$.

Similarly, the external differential
$d:\Omega^*(M)\to\Omega^{*+1}(M)$, gives rise to a
map $\tid:\Sigma\Omega(\mc)\to\Sigma\Omega(\mc),\ \tid s\alpha=-sd\alpha$.
The corresponding coderivative on
$\Lambda^*\Sigma\Omega(\mc)$ is denoted by $D$, its bi-degree is $(0,1)$.

\begin{prop}
\label{propdli}
The maps $L_X,\ I_X$ and $D$ verify the same Cartan identities as the
original maps on the level of $\Omega(\mc)$, i.e.
\begin{align*}
[L_X,\,L_Y]&=L_{[X,\,Y]},&  [I_X,\,I_Y]&=0,&  [L_X,\,I_Y]&=I_{[X,\,Y]}, \\
  [D,\,L_X]&=0,&  [D,\,I_X]&=L_X,&  {\rm and}\ D^2&=0.
\end{align*}
As above the symbol $[\,,\,]$ denotes the bi-graded
commutator of the coderivatives. Observe, that similar formulae hold for the modified
commutator $\{,\,\}$ (this is due to its definition).
\end{prop}
\begin{proof}
Direct calculations with the help of the formula \eqref{eqnew} and \eqref{eqF},
which in the case of the maps $f:\Lambda^1V\to\Lambda^1V$ reduces to
\begin{equation*}
\begin{split}
\tilde f(a_1\wedge&\dots\wedge a_n)\\
                                  &=\sum_{i=1}^n(-1)^{|a_i|(|a_1|+\dots+|a_{i-1}|)+i-1}f(a_i)\wedge\dots\wedge\widehat a_i\wedge\dots\wedge a_n\\
                                  &=\sum_{i=1}^n(-1)^{|f|(|a_1|+\dots+|a_{i-1}|)}a_1\wedge\dots\wedge f(a_i)\wedge\dots\wedge a_n
\end{split}
\end{equation*}
\end{proof}
\begin{rem}\rm
\label{remdeg}
As a matter of fact, due to the discussion preceding Definition \ref{defL},
one can replace all the considered maps by their ``skewed versions" and
from equation \eqref{eqsigns} it follows, that these ``skewed" maps verify
equations, similar to the equations of Proposition \ref{propdli}, but with
the commutator $[\,,\,]_{tot}$ substituted for $[\,,\,]$. The same remark is
true for all the examples, that will follow. In effect, one could rewrite
these examples and the rest of the paper with ``skewed" operators substituted
for the usual ones and with $[\,,\,]_{tot}$ instead of $[\,,\,]$.
\end{rem}
\paragraph{\fbox{Multiplication}}
Let $\cdot:\Omega(\mc)\otimes\Omega(\mc)\to\Omega(\mc),\ \alpha\otimes\beta\to\alpha\cdot\beta$ be the
external multiplication of forms. It turns out, that it can be extended to
a correctly defined map $m:\Lambda^2\Sigma\Omega(\mc)\to\Lambda^1\Sigma\Omega(\mc)$. Namely, put
$$
m(s\alpha\wedge s\beta)=(-1)^{|\alpha|}s(\alpha\cdot\beta)\quad
{\rm for\ all}\ s\alpha,s\beta\in\Sigma\Omega(\mc).
$$
Observe, that here $|\alpha|$
denotes the {\em non-suspended\/} degree of $\alpha$. One easily
checks, that this definition is compatible with the commutator
relations in $\Lambda^2\Sigma\Omega(\mc)$. Note, that the second degree of
this map is equal to $1$. So, we obtain a homogeneous coderivative
$M$ of $\Lambda^*\Sigma\Omega(\mc)$ of bi-degree $(-1,1)$.
\begin{prop}
\label{crprop}
The map $M$ commutes with
$I_x,\ L_X$ and $D$ (the commutator is understood
in the bi-graded sense, i.e. all the bi-graded commutators, defined as above, vanish).
\end{prop}
\begin{proof}
We shall check the statement solely for the maps $I_X$. In the case of the maps $L_X$ and $D$ the proof is similar.
In the view of the formula \eqref{eqnew}, it is enough to check that the
maps $m$ and $\tim_X$ anti-commute (the sign $\epsilon$ in this case is
equal to $1$).

We compute for all $s\alpha,\ s\beta \in\Sigma\Omega(\mc)$:
\begin{equation}
\label{newereq}
\begin{split}
\tilde\imath_X(m(s\alpha\wedge s\beta))&=(-1)^{|\alpha|}\tilde\imath_X(s(\alpha\cdot\beta))\\
                                              &=(-1)^{|\alpha|+1}s\imath_X(\alpha\cdot\beta)\\
                                              &=(-1)^{|\alpha|+1}s\bigl(\imath_X(\alpha)\cdot\beta+(-1)^{|\alpha|}\alpha\cdot\imath_X(\beta)\bigr)
\end{split}
\end{equation}
But $s(\imath_X(\alpha)\cdot\beta)=(-1)^{|\alpha|+1}m(s\imath_X(\alpha)\wedge s\beta)=(-1)^{|\alpha|}m(\tim_X(s\alpha)\wedge s\beta)$.
Similarly $s(\alpha\cdot\imath_X(\beta))=(-1)^{|\alpha|} m(s\alpha\wedge s\imath_X(\beta))=(-1)^{|\alpha|+1} m(s\alpha\wedge\tim_X(s\beta))$
Hence, we continue the equation \eqref{newereq}:
\begin{equation*}
\begin{split}
\tim_X(m(s\alpha\wedge s\beta))&=-m(\tim_X(s\alpha)\wedge
s\beta)+(-1)^{|\alpha|}m(s\alpha\wedge\tim_X(s\beta))\\
                               &=-(m(\tim_X(s\alpha)\wedge s\beta)+(-1)^{(|\alpha|-1)(|\beta|-1)+1}m(\tim_X(s\beta)\wedge
s\alpha)).
\end{split}
\end{equation*}
\end{proof}
\begin{rem}\rm
\label{rem}
As a matter of fact, in the proof of the proposition \label{crprop} we used
only the fact that $\imath_X$ is a differentiation of the algebra $\Omega(\mc)$.
It is easy to check that this proposition holds for arbitrary
differentiations of the de~Rham algebra of a manifold. Namely: let
$\phi:\Omega(\mc)\to\Omega(\mc)$ be a degree $k$ differentiation, i.e.
$\phi(\Omega^i(\mc))\subseteq\Omega^{i+k}(\mc)$, and
$$
\phi(\alpha\cdot\beta)=\phi(\alpha)\cdot\beta+(-1)^{k|\alpha|}\alpha\cdot\phi(\beta).
$$
Define a degree-$k$ map $\tilde\phi:\Sigma\Omega(\mc)\to\Sigma\Omega(\mc)$
as $\tilde\phi(s\alpha)=(-1)^ks\phi(\alpha)$. Then the coderivative $\Phi$, extending this map to
$\Lambda\Sigma\Omega(\mc)$, commutes with the map $M$ above in the bi-graded sense.
\end{rem}
\paragraph{\fbox{The cup-product and its properties}}
The following definition is important for our proof of the formality.
\begin{df}
Let $\phi:\Lambda^k\Sigma\Omega(\mc)\to\Sigma\Omega(\mc)$ and
$\psi:\Lambda^l\Sigma\Omega(\mc)\to\Sigma\Omega(\mc)$ be two linear maps. We
define map $\phi\cup\psi:\Lambda^{k+l}\Sigma\Omega(\mc)\to\Sigma\Omega(\mc)$
by the following formula (here $a_1,\dots a_{k+l}$ are arbitrary elements of $\Sigma\Omega(\mc)$)
\begin{equation}
\begin{split}
&\phi\cup\psi(a_1\wedge\dots\wedge
a_{k+l})\\
&\ =\frac{1}{k!l!}\sum_{\sigma\in\Sigma_{k+l}}(-1)^{\sigma'}\phi(a_{\sigma(1)}\wedge\dots\wedge
a_{\sigma(k)})\cdot\psi(a_{\sigma(k+1)}\wedge\dots\wedge a_{\sigma(k+l)})
\end{split}
\end{equation}
where the sign
$(-1)^{\sigma'}=(-1)^{\sigma+(|\psi|+1)(\sum_{i=1}^k|a_{\sigma(i)}|)+(k-1)(l-1)+|\phi|}$
(the sign $(-1)^\sigma$ is defined above).
\end{df}

\begin{lemma}\ \\
\vspace{-.7cm}
\begin{description}
  \item[{\rm ({\it i\/})}] The cup-product is (bi)graded-commutative, i.e.
\begin{equation}
\phi\cup\psi=(-1)^{|\phi||\psi|+(k-1)(l-1)}\psi\cup\phi.
\end{equation}
  \item[{\rm ({\it ii\/})}] For any differentiation $\delta:\Omega(\mc)\to\Omega(\mc)$, the following formula
holds
\begin{equation}
[\tilde\delta,\,\phi\cup\psi]=[\tilde\delta,\,\phi]\cup\psi+(-1)^{|\delta||\phi|}\phi\cup[\tilde\delta,\,\psi],
\end{equation}
where $\tilde\delta:\Sigma\Omega(\mc)\to\somc$ is defined in remark \ref{rem}
and commutator $[\,,\,]$ -- by formula \eqref{eqnew}.
  \item[{\rm ({\it iii\/})}] If $\alpha,\ \beta,\ \gamma$ and $\delta$ are differentiations of
$\Omega(\mc)$, then
\begin{equation}\label{eqqqq}
\begin{split}
[\tilde\alpha\cup\tilde\beta,\,\tilde\gamma\cup\tilde\delta]&=\tilde\alpha\cup[\tilde\beta,\,\tilde\gamma\cup\tilde\delta]+(-1)^{|\beta|(|\gamma|+|\delta|)}[\tilde\alpha,\,\tilde\gamma\cup\tilde\delta]\cup\tilde\beta\\
&=\ta\cup[\tb,\,\tg]\cup\td+(-1)^{|\alpha||\gamma|}\ta\cup\tg\cup[\tb,\,\td]\\
&\quad+(-1)^{|\beta|(|\gamma|+|\delta|)+|\alpha||\gamma|}\tg\cup[\ta,\,\td]\cup\tb\\
&\quad+(-1)^{|\beta|(|\gamma|+|\delta|)}[\ta,\,\tg]\cup\td\cup\tb.
\end{split}
\end{equation}
\end{description}
\end{lemma}
\begin{proof}
Parts ({\it i\/}) and ({\it ii\/}) follow by direct inspection of formulae.
Part ({\it iii\/}) follows from ({\it ii\/}) and the fact that the map
$\ta\cup\tb(sa\wedge-):\somc\to\somc$ (where $-$ stands for the argument and $sa\in\somc$ is an arbitrary
element) is generated by a degree $|\alpha|+|\beta|+|a|+1$ differentiation
\begin{equation*}
\begin{split}
\imath_a(\alpha\cup\beta)(x)&\eqdef(-1)^{(|\beta|+1)(|a|+1)+|\alpha|}\alpha(a)\cdot\beta(x)\\
&\quad+(-1)^{(|a|+1)(|x|+1)+(|\beta|+1)(|x|+1)+|\alpha|+1}\alpha(x)\cdot\beta(a)
\end{split}
\end{equation*}
of $\Omega(\mc)$.
\end{proof}
\subsection{Proof of the main theorem}
Let $\pi=\sum_kX_k\wedge Y_k$ be a bivector on $\mc$. One easily checks
that the map $\tilde{\tilde\pi}=\sum_k\tim_{X_k}\cup\tim_{Y_k}$ is well-defined, i.e. doesn't
depend on the choice of $X_k,\ Y_k$ in the presentation of $\pi$. The
following statement is evident.
\begin{lemma}\label{lemmm}\ \\
\vspace{-.7cm}
\begin{description}
  \item[{\rm({\it i\/})}] Let $sa,\ sb\in\somc$ be arbitrary elements, then
$$
\tilde{\tilde\pi}(sa\wedge sb)=(-1)^{|a|}s\tilde\pi(a,\,b),
$$
where $\tilde\pi$ is the map defined in equation \eqref{eq5 3/4}.
  \item[{\rm({\it ii\/})}] Let us denote the coderivation of
$\Lambda\somc$, induced by the map $\tilde{\tilde\pi}=\tilde\pi$
by $\Pi$. The following formula holds
$$
[\tilde{\tilde\pi},\,\tid](sa\wedge sb)=s\{a,\,b\},
$$
where $\{\,\,,\,\}$ is the Poisson bracket from the first paragraph.
\end{description}
\end{lemma}
\begin{cor}\label{corr}
The map $s\{\}$, defined as $sa\wedge sb\stackrel{s\{\}}{\mapsto} s\{a,\,b\}$
is equal to
$$s\{\}=\sum_k\bigl(\tim_{X_k}\cup\tlc_{Y_k}-\tlc_{X_k}\cup\tim_{Y_k}\bigr).$$
\end{cor}
 Let us denote by $e^\Pi$ the map
\begin{equation}
1+\Pi+\frac12\Pi\circ\Pi+\frac16\Pi\circ\Pi\circ\Pi+
\ldots:\Lambda\somc\to\Lambda\somc
\end{equation}
It is easy to check that this map is a homomorphism of coalgebras
(in effect this follows just from the fact, that $\Pi$ is a coderivative).
The following formula is the principal result of this paper.
\begin{Th}
If the bivector $\pi$ is a Poisson bivector (i.e.
its Schouten bracket with itself vanishes), then
\begin{equation}
\label{eee}
e^\Pi\circ D\circ e^{\Pi}=D+\widetilde{\{\,\,,\,\}},
\end{equation}
where $\widetilde{\{\,\,,\,\}}$ is the extension to $\Lambda\somc$ of the
Poisson bracket.
\end{Th}
\begin{proof}
Consider the map $e^{t\Pi}$, where $t$ is the formal parameter, $t\in\mathbb
R$ and the formal deformation $e^{t\Pi}\circ D\circ e^{t\Pi}$ of the left
hand side of the formula \eqref{eee}. We shall prove, that it coincides with
the formal deformation $D+t\widetilde{\{\,\,,\,\}}$ of the right hand side of the same
formula, i.e. we are going to prove the formula
$$
e^{t\Pi}\circ D\circ e^{t\Pi}=D+t\widetilde{\{\,\,,\,\}}\eqno{(\ref{eee}')}
$$
To this end it is enough to show, that the formal derivatives at
zero $\frac{\partial^n}{\partial t^n}|_{t=0}$ of the deformed maps on
both sides are equal for all $n$.

Indeed, the equality is evident when $n=0$. If $n=1$, we compute (we omit the
composition signs, where it is possible):
\begin{multline*}
\frac{\partial}{\partial t}(e^{t\Pi}\circ D\circ e^{t\Pi})\\
=e^{t\Pi}(\Pi\circ D+D\circ\Pi)e^{t\Pi}=e^{t\Pi}[\Pi,\,D]e^{t\Pi}=e^{t\Pi}\widetilde{[\tilde{\tilde\pi},\,\tid]}e^{t\Pi}=e^{t\Pi}\widetilde{\{\,\,,\,\}}e^{t\Pi}.
\end{multline*}
We have used the Lemma \ref{lemmm} and the fact that in our case, the bidegrees of the maps being $(1,0)$ and $(-2,1)$, the bi-graded commutator coincides with the anti-commutator (the bi-degree of $D$ is equal to $(1,0)$ and that of $\Pi$ -- to $(-1,-1)$). So, $\frac{\partial}{\partial t}|_{t=0}$
applied to both sides of the formula (\ref{eee}$'$) gives
$\widetilde{\{\,\,,\,\}}$.

For $n=2$ similar calculations give
$$
\frac{\partial^2}{\partial t^2}(e^{t\Pi}\circ D\circ e^{t\Pi})
=\frac{\partial}{\partial
t}(e^{t\Pi}\widetilde{\{\,\,,\,\}}e^{t\Pi})=e^{t\Pi}\widetilde{[\tilde{\tilde\pi},\,s\{\}]}e^{t\Pi}.
$$
Now we can use the formula \eqref{eqqqq} and the Cartan relations
verified by the Lie derivatives and the internal derivatives to show, that
(up to a certain sign) the commutator $[\tilde{\tilde\pi},\,s\{\}]$ is equal
to
\begin{equation*}
\begin{split}
&\sum_{k,j}\{\tim_{X_k}\cup\tim_{[Y_k,\,X_j]}\cup\tim_{Y_j}-\tim_{X_k}\cup\tim_{X_j}\cup\tim_{[Y_k,\,Y_j]}\\
&\quad-\tim_{[X_k,\,X_j]}\cup\tim_{Y_k}\cup\tim_{Y_j}+\tim_{X_j}\cup\tim_{[X_k,\,Y_j]}\cup\tim_{Y_k}\}
\end{split}
\end{equation*}
(of course, some of the terms may cancell each other).
This expression depends only on the Schouten bracket of $\pi$ with itself. Actually, the
Schouten bracket $\{\pi,\,\pi\}$ is defined with the help of the commutation
relations, similar to \eqref{eqqqq}, and one checks that it is equal to
\begin{equation*}
\begin{split}
&\sum_{k,j}\{X_k\wedge[Y_k,\,X_j]\wedge Y_j-X_k\wedge X_j\wedge[Y_k,\,Y_j]\\
&\quad-[X_k,\,X_j]\wedge Y_k\wedge Y_j+X_j\wedge[X_k,\,Y_j]\wedge Y_k\},
\end{split}
\end{equation*}
 and the map, given by the formula
$\sum_{m}\tim_{X_m}\cup\tim_{Y_m}\cup\tim_{Z_m}$ depends only the
tri-vector $\sum_m X_m\wedge Y_m\wedge Z_m$.

From these observations and from the hypotheses that the Schouten bracket of $\pi$ with itself
is equal to $0$, it follows, that $\frac{\partial^2}{\partial t^2}(e^{t\Pi}\circ D\circ e^{t\Pi})=0$.
Consequently, one has $\frac{\partial^3}{\partial
t^3}(e^{t\Pi}\circ D\circ e^{t\Pi})=\ldots=\frac{\partial^n}{\partial
t^n}(e^{t\Pi}\circ D\circ e^{t\Pi})=\dots=0$.
\end{proof}
\begin{rem}\rm
It would, probably be more convenient and more insightful to use the ``skewed" version of all
the operators. Then the bi-graded commutator $[\,,\,]$ should be replaced by
$[\,,\,]_{tot}$ (see remark \ref{remdeg}). Consequently, formula \eqref{eee}
would look like $e^{\bar\Pi}\circ\bar D\circ e^{-\bar\Pi}=\bar
D+\overline{\widetilde{\{\,\,,\,\}}}$, i.e. the operator $e^{\bar\Pi}$
intertwines the standard $L_\infty$\/-structure on $\Sigma\Omega(\mc)$ (that
is the structure, induced by the DGL-algebra structure on it) with the
trivial one (that is the structure, for which all the maps $l_i$ are equal
to $0$, when $i\ge2$).
\end{rem}
\paragraph{\fbox{Corrollaries and discussions}}
\begin{cor}
The Lie algebra $\Omega(\mc)$ is formal.
\end{cor}
\begin{proof}
We need to produce an $L_\infty$\/-quasi-isomorphism between
$\Sigma\Omega(\mc)$ and $\Sigma H(\mc)$. Since all such maps are invertible,
it is enough to find a morphism only in one direction. So, we shall construct
an $L_\infty$\/-quasi-isomorphism from  $\Sigma H(\mc)$ to
$\Sigma\Omega(\mc)$.

Let us first rewrite \eqref{eq10}, taking in account that both the differential and
the Lie bracket vanish in${\Sigma H(\mc)}$:
\begin{equation}
\label{eq11}
\begin{split}
df_n(&\a_1\w\dots\w\a_n)=\\
&=\sum_{i\ge[\frac{n}{2}],\sigma\in
Sh_{i,n-i}}(-1)^{\sigma+(i-1)\varepsilon_i}\{f_i(\a_{\sigma(1)}\w\dots\w\a_{\sigma(i)}),f_{n-i}(\a_{\sigma(i+1)}\w\dots\w\a_{\sigma(n)})\}.
\end{split}
\end{equation}
for all $\a_i\in {\Sigma H(\mc)}$. Or, in brief $(D+\widetilde{\{\,\,,\,\}})\circ F=0$,
where $F$ is the map of coalgebras $F:\wedge^*\Sigma
H(\mc)\to\wedge^*\Sigma\Omega(\mc)$, assembled from $f_k$. It is the map
$F$ that we shall construct.

To this end we first define the bottom stage of this morphism as a linear
splitting of the projection $Z(\mc)\to H(\mc)$ from the space of closed forms to cohomology (we omit the suspension signs):
$f_1([\a])=\a$, where $\a$ is a closed form, representing
the class $[\a]$. For instance, one can choose $\a$ to be the only
harmonic form (with respect to a Riemannian structure) in the class $[\a]$.
We shall have $df_1([\a])=0$, as prescribed by \eqref{eq11}, and it is
almost by definition, that $f_1$ is a quasi-isomorphism of chain complexes.
But the commutator of two harmonic forms is not necessarily harmonic, so
$f_1$ is not in general a homomorphism.

Now we can extend $f_1$ in a trivial way to the map of coalgebras $\wedge^*f_1:\wedge^*\Sigma\Omega(\mc)\to\wedge^*
H(\mc)$. It is clear that $D\circ\wedge^*f_1=0$.

Let us put $F=e^{-\Pi}\circ\wedge^*f_1$. Then
$$
(D+\widetilde{\{\,\,,\,\}})\circ F= (D+\widetilde{\{\,\,,\,\}})\circ
e^{-\Pi}\circ\wedge^*f_1=e^{\Pi}\circ D\circ\wedge^*f_1=0.
$$
\end{proof}

One can visualize the first few stages of $F$. For instance,
\begin{align}
\label{forex1}
f_2([\a]\w\,[\b])&=\widetilde\pi(f_1([\a]),\,f_1([\b])),\\
\label{forex2}
f_3([\a_1]\w[\a_2]\w[\a_3])&=\sum_{\sigma\in\Sigma_3}\frac12\widetilde\pi(\widetilde\pi(f_1([\a_{\sigma(1)}]),f_1([\a_{\sigma(2)}])),f_1([\a_{\sigma(3)}])).
\end{align}

In effect, these formulas can be found independently. For instance, it
follows from \eqref{eq5 1/2} that
$$
\{f_1([\a]),\,f_1([\b])\}=\{\a,\,\b\}=d\widetilde\pi(\a,\,\b).
$$
(Here $\a,\ \b$ denote the forms, representing classes $[\a],\ [\b]$.)
So we have
$$
\{f_1([\a]),\,f_1([\b])\}=df_2([\a]\wedge\,[\b]),
$$
as prescribed by \eqref{eq11}.

Finally for this choice of $f_1,\ f_2$ and $n=3$ \eqref{eq11} takes the following
form:
\begin{equation}
\label{eq12}
\begin{split}
df_3([\a]\w[\b]\w[\g])&=\{f_2([\a]\w[\b]),f_1([\g])\}+(-1)^{|\a|(|\b|+|\g|)}\{f_2([\b]\w[\g]),f_1([\a])\}\\
                      &\qquad+(-1)^{|\b||\g|}\{f_2([\a]\w[\g]),f_1([\b])\}\\
                      &=\{\widetilde\pi(\a,\b),\g\}+(-1)^{|\a|(|\b|+|\g|)}\{\widetilde\pi(\b,\g),\a\}\\
                      &\qquad+(-1)^{|\b||\g|}\{\widetilde\pi(\a,\g),\b\}
\end{split}
\end{equation}
(here $\a,\ \b,\ \g$ are some closed forms).
Let us introduce a map $\pi_3$ by the following equation (for all forms $\a,\ \b,\
\g$):
\begin{equation}
\begin{split}
\pi_3(\a,\b,\g)&=\{\widetilde\pi(\a,\b),\g\}+(-1)^{|\a|(|\b|+|\g|)}\{\widetilde\pi(\b,\g),\a\}\\
                      &\qquad+(-1)^{|\b||\g|}\{\widetilde\pi(\a,\g),\b\}.
\end{split}
\end{equation}
This map extends to all $\Sigma\Omega(\mc)$ the right hand side of
\eqref{eq12}. One can show, that $\pi_3$ verifies the following equation:
\begin{equation}
\label{eq1000}
\begin{split}
\pi_3(\a_1\w\a_2,\b,\g)&=\a_1\w\pi_3(\a_2,\b,\g)\pm\pi_3(\a_1,\b,\g)\w\a_2\\
                       &\quad\pm\Bigl(\{\a_1,\,\b\}\w\widetilde\pi(\a_2,\,\g)+(-1)^{|\a_1|+|\b|}\widetilde\pi(\a_1,\,\b)\w\{\a_2,\,\g\}\Bigr)\\
                       &\quad\pm\Bigl(\{\a_1,\,\g\}\w\widetilde\pi(\a_2,\,\b)+(-1)^{|\a_1|+|\g|}\widetilde\pi(\a_1,\,\g)\w\{\a_2,\,\b\}\Bigr).
\end{split}
\end{equation}
Observe, that if all the forms $\a_1,\ \a_2,\ \b,\ \g$ are closed
then it follows from \eqref{eq5 1/2} that the last two lines on the right hand side of
\eqref{eq1000} consist of exact forms. Now, one can consider local coordinates and
use the argument, similar to the first proof of \eqref{eq5 1/2} to show,
that the following is true on the whole manifold
$$
\pi_3(\a_1,\a_2,\a_3)=d\sum_{\sigma\in\Sigma_3}\frac12\widetilde\pi(\widetilde\pi(\a_{\sigma(1)},\a_{\sigma(2)}),\a_{\sigma(3)}),
$$
which means that \eqref{forex2} is a good choice of $f_3$.
One can continue this process infinitely to find all the maps $f_i,\
i\ge4$ in a straightforward way.

\end{document}